\documentclass[12pt]{article}
\usepackage{amssymb, amsfonts, amsmath, amscd, theorem}
\usepackage[all]{xy}

\numberwithin{equation}{section}
\newcommand{\Ind}{\operatorname{Ind}}
\newcommand{\res}{\operatorname{res}}
\newcommand{\Ker}{\operatorname{Ker}}

\newcommand{\Hom}{\operatorname{Hom}}
\newcommand{\End}{\operatorname{End}}
\newcommand{\Ext}{\operatorname{Ext}}
\newcommand{\Ann}{\operatorname{Ann}}

\newcommand{\gr}{\operatorname{gr}}
\newcommand{\grmod}{\operatorname{grmod}}
\newcommand{\tto}{\twoheadrightarrow}

\newcommand{\mood}{\hbox{\ensuremath{
\operatorname{mod}}}}

\newcommand{\ad}{\hbox{\ensuremath{\operatorname{ad}}}}

\newcommand{\Proj}{\ensuremath{\mathbf{Proj}}}
\newcommand{\Z}{\ensuremath{\mathcal{Z}}}
\newcommand{\F}{\ensuremath{\mathcal{F}}}
\newcommand{\C}{\ensuremath{{k}}}
\newcommand{\D}{\ensuremath{\mathcal{D}_q}}
\newcommand{\Dl}{\ensuremath{\mathcal{D}^\lambda_q}}

\newcommand{\N}{\ensuremath{\mathbb{N}}}

\newcommand{\g}{\ensuremath{\mathfrak{g}}}

\newcommand{\h}{\ensuremath{\mathfrak{h}}}

\newcommand{\bb}{\ensuremath{\mathfrak{b}}}

\newcommand{\BGG}{\ensuremath{\mathcal{O}}}
\newcommand{\Oq}{\ensuremath{\mathcal{O}_q}}
\newcommand{\Oloc}{\ensuremath{\mathcal{O}_{q,loc}}}

\newcommand{\M}{\ensuremath{\mathcal{M}}}
\newcommand{\MBG}{\ensuremath{\mathcal{M}_{B_q}(G_q)}}
\newcommand{\MBpt}{\ensuremath{\mathcal{M}_{B_q}}}

\newcommand{\DBG}{\ensuremath{\mathcal{D}^\lambda_{B_q}(G_q)}}
\newcommand{\DOBG}{\ensuremath{\mathcal{D}^0_{B_q}(G_q)}}

\newcommand{\U}{\ensuremath{U^{fin}_q}}
\newcommand{\Ul}{\ensuremath{U^{\lambda}_q}}
\newcommand{\Uzero}{\ensuremath{U^{0}_q}}
\newcommand{\Dzero}{\ensuremath{\mathcal{D}^0_q}}

\theoremstyle{plain}

\newtheorem{Thm}{Theorem}[section]
\newtheorem{Prop}[Thm]{Proposition}
\newtheorem{Lem}[Thm]{Lemma}

\newtheorem{Cor}[Thm]{Corollary}

\theoremstyle{definition}
\newtheorem{defi}[Thm]{Definition}

\theoremstyle{remark}
\newtheorem{Rem}[Thm]{Remark}

\begin{document}
\title{Quantum flag varieties, equivariant quantum $\mathcal{D}$-modules,  and localization of Quantum groups.}

\author{Erik Backelin and Kobi Kremnitzer}
\maketitle

\begin{abstract} Let $\Oq(G)$ be the algebra of quantized
functions on an algebraic group $G$ and $\Oq(B)$ its quotient
algebra corresponding to a Borel subgroup $B$ of $G$. We define
the category of sheaves on the "quantum flag variety of $G$" to be
the $\Oq(B)$-equivariant $\Oq(G)$-modules and prove that this is a
proj-category. We construct a category of equivariant quantum
$\mathcal{D}$-modules on this quantized flag variety and prove the
Beilinson-Bernsteins localization theorem for this category in the
case when $q$ is not a root of unity.
\end{abstract}

\section{Introduction}  Let $k$ be a field of characteristic zero and fix $q \in k^\star$.
Let $\g$ be a semi-simple Lie algebra over $k$ and let $G$ be the
corresponding simply connected algebraic group. Let $U_q$ be a
quantized enveloping algebra of $\g$. Let $\Oq$ be the algebra of
quantized functions on $G$. Let $\Oq(B)$ be the quotient Hopf
algebra of $\Oq$ corresponding to a Borel subgroup $B$ of $G$.

Adopting Grothendieck's philosophy that a space is the same thing
as its category of sheaves, we define the "quantized flag variety
of $G$", denoted $\MBG$, to be the category of
$\Oq(B)$-equivariant $\Oq$-modules. Thus, an object of $\MBG$ is a
left $\Oq$-module $M$ equipped with a right $\Oq(B)$-coaction such
that the action map is a morphism of $\Oq(B)$-comodules, see
definition \ref{hejsan}. In this language, the global section
functor $\Gamma: \MBG \to \C-\mood$ is the functor of taking
$\Oq(B)$-coinvariants.

Due to Serre's theorem a projective variety can be described
completely algebraically as the quotient of the category of graded
modules over a graded ring modulo its subcategory of torsion
modules. In particular, the category $\M(G/B)$ of quasi-coherent
sheaves on the flag variety $G/B$ is isomorphic to the category
$\Proj(\BGG(G/N))$, where $\BGG(G/N)$ is the algebra of functions
on the basic affine space $G/N$, $N$ is the unipotent radical of
$B$.

A main idea in the theory of non-commutative geometry, due to
Gabriel, Artin and Zhang and others is that this construction
generalizes to non-commutative algebras. The algebra $\BGG(G/N)$
is the so called representation ring of $\g$ and quantizes
naturally to an algebra $\BGG_q(G/N)$ . Lunts and Rosenberg who
were the first to study quantized rings of differential operators
on flag varieties takes $\Proj(\BGG_q(G/N))$ as a definition for
the category of quantized sheaves on $G/B$. We prove in
proposition \ref{aff} that our definition is equivalent to theirs.
The essential thing to prove is proposition \ref{aff}, which
states that $\Oq(\rho)$ is ample. This is not difficult, but much
more complicated than the classical case where one simply uses an
embedding of $G/B$ into a suitable $\mathbf{P}^n$ (there are many
different quantized $\mathbf{P}^n$ and they are not easy to deal
with for this purpose). The key ingredients in our proof is
Kempf-vanishing of Andersen, Polo and Wen Kexin \cite{APW} and a
quantized version of the fact that a space $G/N$ is quasi affine
if and only if every rational $N$-module embeds $N$-linearly to a
rational $G$-module.

Once this technical difficulty is overcome it turns out that our
equivariant sheaves are much easier to deal with than the
proj-approach. Actually, except for section \ref{commutative}
which concerns bimodule structures on $\Oq(G)$-equivariant sheaves
and is independent of the bulk material of this paper, we don't
have any explicit need of the proj-category, but we frequently use
the fact that $\Oq(\rho)$ is ample.
\smallskip

\noindent In particular this becomes evident in the study of
$D$-modules: It is not clear what a quantized ring of differential
operators should be on a non-commutative ring (and even less so on
a non-commutative space).  Lunts and Rosenberg \cite{RL1} gave a
definition of such a ring of differential operators, using a
definition similar to Grothendieck's classical construction, that
works for any graded algebra once they fixed a certain
bi-character on it. This construction has the disadvantage that
ring of differential operators it produces seems to be too big.
They apply this construction to $\Oq(G/N)$ (see \cite{RL2}) and
define a $D$-module on (quantum) $G/B$ to be an object in the
quotient category of graded $D$-modules on $G/N$ modulo torsion
modules.

Recently, Tanisaki \cite{T} defined the ring of differential
operators on quantum $G/N$ to be the subalgebra of
$\End_k(\Oq(G/N))$ generated by $\Oq(G/N)$ and $U_q$. This is a
subalgebra of Lunts and Rosenberg's algebra of differential
operators. The category of $\mathcal{D}$-modules he gets on
quantum $G/B$ by the proj-construction is equivalent to the one we
get.

In our equivariant approach we don't need a ring of differential
operators on $G/N$; we only need a ring $\D$ of differential
operators on $G$ and we simply define $\D$ to be the smash product
algebra $\Oq \star U_q$. We define a $\lambda$-twisted quantum
$D$-module on $G/B$ ($\lambda$ is an element in the character
group of the weight lattice of $\g$) to be an object $M \in \MBG$
with an additional action of $\D$ such that the coaction of
$\Oq(B)$ and the action of $U_q(\bb) \subset U_q \subset \D$ on
$M$ "differs by $\lambda$". The $\lambda$-twisted $D$-modules
forms a category denoted $\DBG$. See definition \ref{d1}.

In the equivariant language, there is no ($\lambda$-twisted) sheaf
of rings of differential operators on $G/B$. But we do have a
distinguished object $\Dl$ which represents the global sections.
It can be described as the maximal quotient of $\D$ that belongs
to $\DBG$. As an object of $\MBG$, $\Dl$ is isomorphic to the
"induced sheaf" $\Oq \otimes M_\lambda$, where $M_\lambda$ is a
Verma module with highest weight $\lambda$. As the global section
functor $\Gamma$ is given by $\Hom(\Dl, \;)$ we see that
$\Gamma(\Dl) = \End(\Dl)$ is an algebra.

\smallskip

\noindent We prove in proposition \ref{globsectD} that for each
$q$ except a finite set of roots of unity (depending on $\g$),
$\Gamma(\Dl)$ is isomorphic to $\U/J_\lambda$, where $J_\lambda$
is the annihilator of $M_\lambda$.

The proof of \ref{globsectD} uses the corresponding classical
result for the case $q = 1$ and results of Joseph and Letzter
\cite{JL2} which states that the standard filtration on the
enveloping algebra $U(\g)$ has a quantized version where the
subquotients have the same dimensions as in the classical case.

\smallskip

\noindent The main result of this paper is theorem \ref{lockie}
which is the quantized version of Beilinson-Bernsteins
localization, \cite{BB}. It states that the global section functor
gives an equivalence between $\DBG$ and the category of modules
over the algebra $\U/J_\lambda$ when $\lambda$ is dominant and
regular. Here $\U$ denotes the ad-finite part of $U_q$. This
theorem holds only if $q$ is not a root of unity and the reason
for this is that Harish Chandra's description of the center of
$U_q$ doesn't hold at a root of unity.  Our proof of the
localization theorem is almost identical to the one given in
\cite{BB}.

\smallskip

\noindent Lunts and Rosenberg \cite{RL2} conjectured proposition
\ref{globsectD} and theorem \ref{lockie} for their $D$-modules and
Tanisaki proved them for his.

\subsection{Acknowledgements} We thank Joseph Bernstein
and David Kazhdan for many useful conversations. The second author
would like to thank the IHES for their hospitality and for
financing his visit the summer of 2003 during which parts of this
work was developed.
\section{Generalities}
\subsection{Quantum groups} See Chari and Pressley \cite{CP} for details
about the topics in this section: Let $k$ be a field of
characteristic zero and fix $q \in k^\star$. Let $\g$ be a
semi-simple Lie
algebra 
and let $\h \subset \bb$ be a Cartan subalgebra contained in a
Borel subalgebra of $\g$. Let $P \subset \h^\star$ be the weight
lattice and $P_+$ the positive weights; the $i$'th fundamental
weight is denoted by $\omega_i$ and $\rho$ denotes the half sum of
the positive roots. Let $Q \subset P$ be the root lattice and $Q_+
\subset Q$ those elements which have non-negative coefficients
with respect to the basis of simple roots. Let $\mathcal{W}$ be
the Weyl group of $\g$. We let $<\,,\,>$ denote a
$\mathcal{W}$-invariant bilinear form on $\h^\star$ normalized by
$<\gamma,\gamma> = 2$ for each short root $\gamma$.

Let $T_P = \Hom_{groups}(P,k^\star)$ be the character group of
$P$. be the character group of $P$ with values in $k$ (we use
additive notation for this group). If $\mu \in P$, then $<\mu,P>
\subset \mathbb{Z}$ and hence we can define $q^\mu \in T_P$ by the
formula $q^\mu(\gamma) = q^{<\mu,\gamma>}$, for $\gamma \in P$. If
$\mu \in P, \lambda \in T_P$ we write $\mu + \lambda = q^\mu +
\lambda$. Note that the Weyl group naturally acts on $T_P$.

Let $U_q$ be the simply connected quantized enveloping algebra of
$\g$ over $k$. Recall that $U_q$ has algebra generators $E_\alpha,
F_\alpha, K_\mu$, $\alpha, \beta$ are simple roots, $\mu \in P$
subject to the relations
$$ K_{\lambda} K_{\mu} = K_{\lambda+\mu},\;\; K_0 = 1,$$
$$K_\mu E_\alpha K_{-\mu} = q^{<\mu, \alpha>} E_\alpha,\;\; K_\mu F_\alpha K_{-\mu} = q^{-<\mu, \alpha>}
F_\alpha,$$
$$[E_\alpha,\, F_\beta] = \delta_{\alpha, \beta} {{K_\alpha - K_{-\alpha}}\over{q-q^{-1}}}$$
and certain Serre-relations that we do not recall here. (We assume
that $q^2 \neq 1$.)

Let $G$ be the simply connected algebraic group with Lie algebra
$\g$, $B$ be a Borel subgroup of $G$ and $N \subset B$ its
unipotent radical. Let $\bb = \operatorname{Lie} B$ and
$\mathfrak{n} = \operatorname{Lie} N$ and denote by $U_q(\bb)$ and
$U_q(\mathfrak{n})$ the corresponding subalgebras of $U_q$. Then
$U_q(\bb)$ is a Hopf algebra, while $U_q(\mathfrak{n})$ is only an
algebra. Let $\Oq = \Oq(G)$ be the algebra of matrix coefficients
of finite dimensional type-1 representations of $U_q$. There is a
natural pairing $(\;,\;): U_q \otimes \Oq \to \C$. This gives a
$U_q$-bimodule structure on $\Oq$ as follows
\begin{equation}\label{23}
ua = a_1 (u,a_2),\;\; au = (u,a_1)a_2, \,\, u \in U_q, a \in \Oq
\end{equation}
Then $\Oq$ is the (restricted) dual of $U_q$ with respect to this
pairing. We let $\Oq(B)$ and $\Oq(N)$ be the quotient algebras of
$\Oq$ corresponding to the subalgebras $U_q(\bb)$ and
$U_q(\mathfrak{n})$ of $U_q$, respectively, by means of this
duality. Then $\Oq(B)$ is a Hopf algebra and $\Oq(N)$ is only an
algebra.

{\it Verma modules:} For each $\lambda \in T_P$ there is the one
dimensional $U_q(\bb)$-module $\C_\lambda$ which is given by
extending $\lambda$ to act by zero on the $E_\alpha$'s. The
Verma-module $M_\lambda$ is the $U_q$-module induced from
$\C_\lambda$. Thus $M_\lambda$ is a cyclic left $U_q$-module with
a generator $1_\lambda$ subject to the relations
\begin{equation}\label{Vermamoduledef}
E \cdot 1_\lambda = 0,\, K_\alpha \cdot 1_\lambda =
\lambda(\alpha) \cdot 1_\lambda.
\end{equation}

Let $\mu \in P$. We write $\C_\mu = \C_{q^\mu}$ and $M_\mu =
M_{q^\mu}$ in this case. Note that $\C_\mu$ integrates to an
$\Oq(B)$-comodule: we can think of $\mu$ as living in the
restricted dual of $U_q(\bb)$ (i.e. in $\Oq(B)$) and $\mu$ is
grouplike. The comodule action on $\C_\mu$ is now given by
\begin{equation}\label{1D}
 1_\mu \to 1_\mu \otimes \mu.
\end{equation}
Each $1$-dimensional $\Oq(B)$-comodule is isomorphic to $\C_\mu$
for some $\mu \in P$.

{\it Harish Chandra homomorphism:} Let $\mathcal{Z}$ denote the
center of $U_q$. Assume that $q$ is not a root of unity. Given
$\lambda \in T_P$ there is the central character $\chi_\lambda:
\mathcal{Z} \to \C$; it is characterized by the property that
$\Ker \chi_\lambda \cdot M_{\lambda - \rho} = 0$. We have $\Ker
\chi_\lambda = \Ker \chi_{w\lambda}$.

Let $\lambda \in T_P$. If $q$ is not a root of unity, we say that
\begin{itemize}
\item{} $\lambda$ is dominant if $\chi_{\lambda} \neq
\chi_{\lambda+\phi}$ for each $\phi \in Q_{+} \setminus \{0\}$.
\item{} $\lambda$ is regular dominant if for all $\phi \in P_+$ and
all weights $\psi$ of $V_\phi$, $\phi \neq \psi$, we have
$\chi_{\lambda +\phi} \neq \chi_{\lambda+\psi}$. (Here $V_\phi$ is
the irreducible finite dimensional type-1 representation of $U_q$
with highest weight $\phi$. See also definition
\ref{representationring}.)
\end{itemize}
If $\lambda = q^\mu, \mu \in P$ this is equivalent to saying that
$\mu$ is dominant, respectively regular dominant, in the usual
sense.

{\it Finite part of $U_q$}. The algebra $U_q$ acts on itself by
the adjoint action $\ad: U_q \to U_q$ where $\ad(u)(v) = u_1 v
S(u_2)$. Let $U^{fin}_q$ be the finite part of $U_q$ with respect
to this action:
$$ U^{fin}_q = \{v \in U_q; \dim \ad(U_q)(v) < \infty\}.$$
This is a subalgebra. (See \cite{JL1}.)

We shall frequently refer to a right (resp. left) $\Oq$-comodule
as a left (resp. right) $G_q$-module, etc. If we have two right
$\Oq$-comodules $V$ and $W$, then $V\otimes W$ carries the
structure of a right $\Oq$-comodule via the formula
$$
\delta(v\otimes w) = v_1\otimes w_1 \otimes v_2w_2$$ We shall
refer to this action as the {\it tensor} or {\it diagonal} action.
A similar formula exist for left comodules.

\subsection{Proj-categories}\label{pro}
We shall use a multigraded version of the classical result about
Proj-categories that is basically due to Serre. We consider tuples
of data $(\mathcal{C};\, $\BGG$;\, s_1, \ldots s_l)$ where
$\mathcal{C}$ is an abelian category, $\BGG$ a fixed object of
$\mathcal{C}$, $s_1, \ldots s_l$ a set of pairwise commuting
autoequivalences of $\mathcal{C}$. For $\mathbf{n} = (n_1, \dots,
n_l) \in \N^l$, and $M \in \operatorname{Ob}(\mathcal{C})$ we
define "twisting-functors" on $\mathcal{C}$ by
$$M(\mathbf{n}) = s^{n_1}_1 \cdots s^{n_l}_l(M).$$

We define for any $M \in \operatorname{Ob}(\mathcal{C})$ its
global sections $\Gamma(M) = \Hom_{\mathcal{C}}(\BGG,F)$. We also
put $\underline{\Gamma}(M) = \oplus_{\mathbf{n} \in \N^l}
\Gamma(M(\mathbf{n}))$.

For any $\Z^l$-graded algebra $R = \oplus_{\mathbf{n} \in \Z^l}
R_{\mathbf{n}}$ we denote by $\Proj(R)$ the quotient category of
the category of $\N^l$-graded left $R$-modules modulo the Serre
subcategory of torsion object. Here, an object is called torsion
if each of its elements is annihilated by $R_{\geq k} =
\oplus_{n_1, \ldots, n_l \geq k} R_{\mathbf{n}}$ for some $k \geq
0$. Let $\mathcal{C}^{0}$ denote the set of noetherian objects in
$\mathcal{C}$. Artin and Zhang \cite{AZ} proved the following
result
\begin{Prop}\label{serre}
Assume that $i)$ $\BGG$ is in $\mathcal{C}^0$;

$ii)$ $\underline{\Gamma}(\BGG)$ is a left-noetherian ring and
$\underline{\Gamma}(M)$ is finitely generated over $\Gamma(\BGG)$
for $M \in \mathcal{C}^0$;

$iii)$ For each $M \in \mathcal{C}^0$ there is a surjection
$\oplus^p_{j=1}\BGG (-\mathbf{n}_j) \to M$; and

$iv)$ if $M,N \in \mathcal{C}^0$ and $M \to N$ is a surjection,
then $\Gamma(M(\mathbf{n})) \to \Gamma(N(\mathbf{n}))$ is
surjective for $\mathbf{n} >> 0$. Then $\mathcal{C}$ is equivalent
to the category $\Proj(\underline{\Gamma}(\BGG))$.
\end{Prop}
We will refer to an autoequivalence satisfying $iii)$ and $iv)$ as
ample.

\section{Quantum flag variety}
\subsection{}
The composition
\begin{equation}\label{e2}
\Oq \to \Oq \otimes \Oq \to \Oq \otimes \Oq(B)\end{equation}
defines a right $\Oq(B)$-comodule structure on  $\Oq$. A
$B_q$-equivariant sheaves on $G_q$ is a triple $(F, \alpha,
\beta)$ where $F$ is a vector space, $\alpha: \Oq \otimes F \to F$
a left $\Oq$-module action and $\beta: F \to F \otimes \Oq(B)$ a
right $\Oq(B)$-comodule action such that $\alpha$ is a right
comodule map, where we consider the tensor comodule structure on
$\Oq(G) \otimes F$.
\begin{defi}\label{hejsan} We denote $\MBG$ to be the category of $B_q$-equivariant
sheaves on $G_q$. Morphisms in $\MBG$ are those compatible with
all structures.\end{defi}

\begin{Rem} In the classical case, when $q = 1$, the category $\M_B(G)$ is
equivalent to the category $\M(G/B)$ of quasi-coherent sheaves on
$G/B$.
\end{Rem}

We similarly have categories $\mathcal{M}(G_q) :=
\mathcal{M}_{\{e\}}(G_q) =$ category of $\Oq$-modules (where
$\{e\}$ is the one-point group) and $\MBpt :=
\mathcal{M}_{B_q}(\operatorname{pt})=$ $B_q$-modules (where
$\operatorname{pt}$ is the one-point space), $\M :=
\M(\operatorname{pt}) = \C-\mood$.

\subsection{}
We have a {\it basic diagram} that will be used throughout this
paper
\begin{equation}\label{basic}
\begin{matrix}
{\mathcal{M}(G_q)} & {\overset{p}{\to}} &
{\M} \\
{{\downarrow}_\pi} & {} & \downarrow_{\bar{\pi}} \\
\MBG & \overset{\bar{p}}{\to} & \MBpt
\end{matrix}
\end{equation}

Here each arrow denotes a pair of adjoint functors; hence the
adjoint pair of functors corresponding to an arrow $f$ will be
denoted $(f^\star,\,f_\star)$ and $f_\star$ goes in the direction
of the arrow.  Here $\pi_\star = (\;)\otimes \Oq(B)$, where $B_q$
acts on the second factor and $\Oq$ acts via the tensor action
(using that $\Oq(B)$ is a quotient of $\Oq$); $\pi^\star =
\operatorname{forget}$; $p_\star = \operatorname{forget}$ and
$p^\star = \Oq\otimes (\;)$, where $\Oq$ acts on the first factor.

Similary, $\overline{\pi}_\star = (\;)\otimes \Oq(B)$, where $B_q$
acts on the second factor; $\overline{\pi}^\star =
\operatorname{forget}$; $p_\star = \operatorname{forget}$;
$\overline{p}^\star  = \Oq \otimes (\;)$ where $\Oq$ acts on the
first factor and $B_q$ acts via the tensor action.

The diagram is commutative in the sense of usual commutativity
after applying lower star (resp. upper star) to all the arrows.
All functors considered are exact; hence all "lower star"
morphisms maps injectives to injectives.

We define
\begin{defi} Let $\lambda \in P$ and put $\Oq(\lambda)
= \overline{p}^\star \C_{-\lambda}$. We call $\Oq(\lambda)$ a line
bundle.
\end{defi}
For each $M \in \MBG$ and $\lambda \in P$ we put
\begin{equation}
M(\lambda) = M \otimes \C_{-\lambda}.
\end{equation} This is an
object in $\MBG$ with the $\Oq$-action on the first factor and the
tensor $B_q$-action called the $\lambda$-twist of $M$.

\subsection{}
\begin{defi} The global section functor $\Gamma: \MBG \to \C-\mood$
is defined by
$$\Gamma(M) = \Hom_{\MBG}(\Oq, M) = \{m \in M; \Delta_B(m) =
m\otimes 1\}.$$ This is the set of $B_q$-invariants in $M$.
\end{defi}

We can now state our main result about the category $\MBG$.

\begin{Prop}\label{aff} $1)$`Each object in $\MBG$ is a quotient of a direct sum of
$\Oq(\lambda)$'s.

$2)$ Any surjection $M \twoheadrightarrow M'$ of noetherian
objects in $\MBG$ induces a surjection $\Gamma(M(\lambda))
\twoheadrightarrow \Gamma(M'(\lambda))$ for $\lambda >> 0$.
\end{Prop}
Here the notation $\lambda >> 0$ means that $<\lambda,
\alpha^\wedge>$ is a sufficiently large integer for each simple
root $\alpha$. Thus, the proposition can be phrased as:
$\Oq(\lambda)$ is ample if $<\lambda, \alpha^\wedge>
> 0$ for each simple root $\alpha$.
\begin{defi}\label{representationring} Let $V_\lambda = \Gamma(\Oq(\lambda))$ and let
$A = \oplus_{\lambda \in P_{+}} V_\lambda$ be the representation
ring of $U_q$.
\end{defi}
Note that the $V_\lambda$'s, $\lambda \in P_+$ are the simple
finite dimensional $U_q$-modules if $q$ is not a root of unity.
\begin{Cor}\label{proj}  The category $\MBG$ is equivalent to
$\Proj (A)$.
\end{Cor}
$\mathit{Proof\; of \; corollary \;\ref{proj}.}$ Then, with the
notations of section \ref{pro}
$$A = \oplus_{\lambda \in P_{+}}\Gamma(\Oq(\lambda))
= \underline{\Gamma}(\Oq).$$
Hence we are left to show that the conditions $i)-iv)$ of
proposition \ref{serre} are satisfied for the tuple $(\MBG, \Oq,
s_1,\ldots , s_l)$, where $s_i(M) = M(\omega_i)$ and we recall
that the $\omega_i$'s are the fundamental weights. Now, $iii)-iv)$
is proposition \ref{aff}; $i)$ holds because $\Oq$ is a noetherian
ring and $ii)$ is clearly true for line bundles and then follows
for general modules from $iii)$. $\Box$
\subsection{}The following two sections are devoted to the proof of proposition \ref{aff}.
Apart from the interesting results corollary \ref{Kempf} and lemma
\ref{GN} the proof consists mostly of rather technical standard
arguments. In this section we show that various categories have
enough injectives and calculate some cohomology groups. We deduce
in corollary \ref{Kempf} that Kempf-vanishing holds in $\MBG$.

Let $M \in \MBG$. The adjunction map $\pi_\star \pi^\star M \to M$
(which is given by $\operatorname{Id} \otimes
\operatorname{counit}$) has a splitting given by the comodule
action $M \to M \otimes \Oq(B) = \pi_\star \pi^\star M$. Let $I$
be an injective hull of $\pi^\star M$ in ${\mathcal{M}(G_q)}$.
Then $M$ embeds into $\pi_\star I$ and we conclude that
\begin{Lem}\label{enough1} The category
$\MBG$ has enough injectives.
\end{Lem}
Let $\tilde{\Gamma}:\MBpt \to \C-\mood$ be the functor of taking
$B_q$-invariants on $\MBpt$. Thus derived functors of $\Gamma$ and
$\tilde{\Gamma}$ are defined. We have

\begin{equation}\label{gamma}
\Gamma = \widetilde{\Gamma} \circ \overline{p}_\star
\end{equation}
The category $\MBpt$ has enough injectives because
$\overline{\pi}_\star$ maps injectives to injectives, each object
in $\mathcal{M}(\operatorname{pt})$ is injective and any $M \in
\MBpt$ imbeds to $\overline{\pi}_\star\overline{\pi}^\star M$. We
have
\begin{Lem}\label{enough2} $1)$ If $I \in \MBpt$ is injective then $\overline{p}^\star I$ is
$\Gamma$-acyclic. $2)$ The functor $\Gamma$ has finite
cohomological dimension and the formula $R\Gamma(M) =
R\widetilde{\Gamma}(\overline{p}_\star M)$ holds for $M \in \MBG$.
\end{Lem}
$\mathbf{Proof}$ $1)$ Let $I \in \MBpt$ be injective. Then $I$
imbeds to $\overline{\pi}_\star \overline{\pi}^\star(I) =
\Oq(B)\otimes I \cong \Oq(B)^{\dim I}$. Since
$\overline{\pi}_\star$ preserves injectives and every object in
$\mathcal{M}$ is injective, $\overline{\pi}_\star
\overline{\pi}^\star(I)$ is injective. Since $I$ is injective this
embedding splits. Thus it suffices to prove that
$\overline{p}^\star (\Oq(B))$ is $\Gamma$-acyclic. We have
$\overline{p}^\star (\Oq(B)) = \pi_\star(\Oq)$ and conclude
$$R^j\Gamma(\overline{p}^\star(\Oq(B))) = \Ext^j_{\MBG}(\Oq, \pi_\star(\Oq))
\cong$$ $$\Ext^j_{\mathcal{M}(G_q)}(\pi^\star(\Oq), \Oq) =
\Ext^j_{\mathcal{M}(G_q)}(\Oq, \Oq),$$  where we used that
$\pi_\star$ is exact and preserves injectives in the second
isomorphism. Since $\Oq$ is projective in $\mathcal{M}(G_q)$ the
last term vanishes for $j
>0$.
\smallskip

\noindent $2)$ Andersen, Polo and Wen Kexin \cite{APW} has shown
that the functor $\Gamma \circ \overline{p}^\star$ has
cohomological dimension $\leq \dim G/B$. Let $M \in \MBG$. Since
$\C$ is a direct summand in $\Oq$, $M$ is a direct summand in
$\overline{p}^\star\overline{p}_\star(M) = \Oq \otimes M$ as a
$B_q$-module. Thus, by $1)$, $R^i\widetilde{\Gamma}(M)$ is a
direct summand in $R^i (\Gamma \circ
\overline{p}^\star)(\overline{p}_\star(M))$ and the latter module
vanishes for $i > \dim G/B$.

Let $M \to I_\bullet$ be an injective resolution in $\MBG$. We get
again
$$
R^i\Gamma(M) = \operatorname{H}^i(\Gamma(I_\bullet)) =
\operatorname{H}^i(\widetilde{\Gamma}(\overline{p}_\star
I_\bullet)) = R^i\tilde{\Gamma}(\overline{p}_\star M)$$ and the
last term vanishes for $i
> \dim G/B$. $\Box$

\begin{Cor}\label{Kempf} [Kempf vanishing] For each $\lambda \in P_{+}$ and each $i > 0$ we have
$R^i\Gamma(\Oq(\lambda) = 0$.
\end{Cor}
$\mathbf{Proof}$ Let $\lambda \in P$. Choose an injective
resolution $\C_\lambda \to I_\bullet$ in $\MBpt$. Then
$$R^i(\Gamma \circ \overline{p}^\star) (\C_\lambda) =
\operatorname{H}^i(\Gamma(\overline{p}^\star I_\bullet)) =
R^i\Gamma(\Oq(\lambda))$$ since the $\overline{p}^\star I_\bullet$
are $\Gamma$-acyclic, by lemma \ref{enough2}. Now, it is shown in
\cite{APW} that $R^i(\Gamma \circ \overline{p}^\star) (\C_\lambda)
= 0$ for $i > 0$, if $\lambda \in P_{+}$. $\Box$

\subsection{} In this section we introduce a $G_q$-equivariant
structure on certain objects in $\MBG$. We prove the key lemma
\ref{GN} and finally we prove proposition \ref{aff}.

Let $V \in \mathcal{M}_{G_q}$. Denote by $V \vert B_q$ the module
$V$ restricted to $B_q$ and by $V^{triv}$ the trivial $B_q$-module
whose underlying space is $V$. We have the following crucial fact
\begin{Lem}\label{l1} The objects $\overline{p}^\star
(V \vert B_q)$ and $\overline{p}^\star(V^{triv})$ are isomorphic
in $\MBG$.
\end{Lem}
{\it Proof.} The map $\overline{p}^\star (V \vert B_q) \to
\overline{p}^\star(V^{triv})$ is given by $a\otimes v \to
av_2\otimes v_1$. It is easily checked that this is an
isomorphism. $\Box$

For any $V \in \mathcal{M}_{B_q}$, $\overline{p}^\star V$ carries
the additional structure of a right $G_q$-module via the (right)
action on the first factor. This structure is compatible with the
left $\Oq$-action and makes $\overline{p}^\star V$ a
$B_q-G_q$-bimodule. We denote by ${}_{G_q}\MBG$ the category of
all objects in $\MBG$ that carry this additional structure. We
have
\begin{Lem}\label{l12} The functor $\overline{p}^\star$ induces an equivalence
$\MBpt \to {}_{G_q}\MBG$.
\end{Lem}
$\mathit{Proof.}$ For $M \in {}_{G_q}\MBG$ denote its
$\Oq$-comodule action by $\Delta$ and let $(M)^{G_q} =\{m \in M;\,
\Delta m = 1 \otimes m\}$ be the set of $G_q$-invariants. It is
straight forward to verify that the functor $(\;)^{G_q}$ is
inverse to $\overline{p}^\star$. $\Box$
\begin{Rem} The map $\overline{p}^\star
(V \vert B_q) \to \overline{p}^\star(V^{triv})$ in lemma \ref{l1}
becomes an isomorphism in ${}_{G_q}\MBG$ if we modify the $G_q$
action on $\overline{p}^\star(V^{triv})$: we define the new
$G_q$-action to be the diagonal action.
\end{Rem}

We first prove the following
\begin{Lem}\label{GN} Assume
$V \in \M_{B_q}$ is finite dimensional and satisfies the
following: if $\C_\lambda$ is a one dimensional sub quotient of
$V$ then $\lambda \in P_{+}$. Then there is a f.d. $W \in
\M_{G_q}$ and an $B_q$-linear surjection $W \twoheadrightarrow V$.
\end{Lem}
$\mathit{Proof \; of \; lemma \; \ref{GN}.}$ We have induction and
restriction functors between categories
\begin{equation}\label{nerd}
\M_{B_q} \underset{\res}{\overset{\Ind}{\rightleftarrows}}
\M_{G_q}
\end{equation}
Let $V \in \M_{B_q}$. We have $\Ind^{G_q}_{B_q}(V) = (\Oq \otimes
V)^{B_q}$. For each one-dimensional sub quotient $\C_\lambda$ of
$V$ the adjunction morphism
$$\res^{B_q}_{G_q}\Ind^{G_q}_{B_q}(\C_\lambda) = \res^{B_q}_{G_q}
(V_\lambda) = V_\lambda \to \C_\lambda$$ is surjective, since
$\lambda \in P_{++}$. An easy induction using corollary
\ref{Kempf} shows that the functor
$res^{B_q}_{G_q}\Ind^{G_q}_{B_q}$ is exact on any sequence of the
form $V' \to V \to V/V'$ for any submodule $V'$ of $V$. By
induction and the five lemma we conclude that
$\res^{B_q}_{G_q}\Ind^{G_q}_{B_q}(V) \to V$ is surjective.

We take $W$ to be any f.d. $G_q$-submodule of
$\Ind^{G_q}_{B_q}(V)$ that surjects to $V$.
\smallskip

\noindent $\mathit{Proof \; of \; proposition \; \ref{aff}.}$ $1)$
Let $M \in \MBG$. We can assume that $M$ is noetherian. Take a
minimal set of generators of $M$ as an $\Oq$-module and let $V$ be
the $B_q$-module they generate; $V$ is f.d. by the noetherian
hypothesis. We get a surjection $\overline{p}^\star V
\twoheadrightarrow M$ in $\MBG$. Take $\lambda \in P$ such that
$V\otimes \C_\lambda$ satisfies the assumption of lemma \ref{GN}
and let $W$ be a f.d. $G_q$-module that surjects to $V\otimes
\C_\lambda$. Then $\Oq \otimes W$ surjects to $(\Oq \otimes
V)(\lambda)$ and hence to $M(\lambda)$. It follows  from lemma
\ref{l1} that $\Oq \otimes W$ is generated by its
$B_q$-invariants. Hence $M(\lambda)$ is as well, i. e. we have a
surjection $\Oq(-\lambda)^m \tto M$.

\smallskip

\noindent $2)$. Let $M \twoheadrightarrow M'$ be a surjection in
$\MBG$. Let $F_0$ be a direct sum of line bundles and $F_0 \tto M$
a surjection. If we can prove that the composition $F_0 \tto M'$
induces a surjection $\Gamma(F_0(\lambda)) \to
\Gamma(M'(\lambda))$ for suitable $\lambda$ it will follow that
the map $\Gamma(M(\lambda)) \to \Gamma(M'(\lambda))$ is surjective
for such $\lambda$ as well.

Put $n = \dim G/B$ which we recall is the cohomological dimension
of the functor $\Gamma$ and pick a resolution
\begin{equation}\label{resolution}
F_n \to \ldots \to F_1 \to F_0 \to M' \to 0
\end{equation}
where each $F_i$ is a direct sum of line bundles. Let $\lambda$ be
sufficiently large for the following property $(\star)$ to hold:
each $F_i(\lambda)$ is a direct sum of various $\Oq(\mu)$, where
each $\mu \in P_+$.  Tensoring \ref{resolution} with
$\C_{-\lambda}$ we get an exact sequence
\begin{equation}\label{resolution2}
F_n(\lambda) \overset{f_n}{\to} \ldots \overset{f_2}{\to}
F_1(\lambda) \overset{f_1}{\to} F_0(\lambda) \overset{f_0}{\to}
M'(\lambda) \overset{f_{-1}} \to 0
\end{equation}
Put $K_i = \Ker f_i$. We must show that $\Gamma(f_0)$ is
surjective. We have short exact sequences $K_i \hookrightarrow
F_i(\lambda) \tto K_{i-1}$ inducing exact sequences
$$
R^i\Gamma(F_{i}(\lambda)) \to R^{i}\Gamma(K_{i-1})  \to
R^{i+1}\Gamma(K_i) \to R^{i+1}\Gamma(F_{i}(\lambda))
$$
By $(\star)$ and corollary \ref{Kempf}, we get isomorphisms
$R^{i}\Gamma(K_{i-1})  \cong R^{i+1}\Gamma(K_i)$, for $i \geq 1$.
Now, $R^{n+1}\Gamma(K_n) = 0$, because $\Gamma$ has cohomological
dimension $n$; hence $R^1\Gamma(K_0) = 0$. Considering the above
sequence when $i = 0$ we conclude that $\Gamma(f_0)$ is
surjective.

\subsection{$G_q$-commutativity of $A$}\label{commutative} The results in this
section are not needed for the rest of this paper.

Classically, a sheaf of $\mathcal{O}_{G/B}$-modules is a bimodule
as $\mathcal{O}_{G/B}$ is commutative. In the quantum case this is
no longer true. Yet the class of $G_q$-equivariant objects in
$\Proj(A)$ admits an $A$-bimodule structure. Using corollary
\ref{proj} one deduces that the $G_q$-equivariant objects in
$\MBG$ act on $\MBG$; we suggestively denote this action by
$\otimes_{\Oq}$.

We recall the notion of a commutative algebra in a braided tensor
category.
\begin{defi} Let $\mathcal{B}$ be a braided tensor category. An
algebra in $\mathcal{B}$ is a pair $(R,m)$ where $R \in
\operatorname{Ob}(\mathcal{B})$ and $m: R\otimes R \to R$
satisfying the usual associativity axiom. $R$ is called
commutative if the diagram
\begin{equation}\label{basic2}
\begin{matrix}
{R\otimes R} & {\overset{m}{\to}} &
{R} \\
{{\downarrow}_\sigma} & {} & \parallel \\
R\otimes R & \overset{m}{\to} & R
\end{matrix}
\end{equation}
commutes, where $\sigma$ is the braiding.
\end{defi}
Similarly, one can define left modules over $R$, etc, in the
braided tensor category. If $R$ is commutative in $\mathcal{B}$
then left modules are bimodules: Let $M$ be a left $R$-module.
Composing the left action with the braiding we get a right action
$$M \otimes R \to R \otimes M \to M$$
It is easily verified that this structure commutes with the left
structure, giving us the asserted bimodule structure.

We now consider the braided tensor category $U_q-\grmod_P$ of
$P$-graded left $U_q$-modules (we assume additionally that each
braided component is finite dimensional). We assume that $q$ has a
square root in $k$ and fix such a root $q^{1/2}$. The braiding in
$U_q-\grmod_P$ is the product of the usual braiding on
$U_q$-modules and the braiding on the category of $P$-graded
vector spaces given by the bicharacter
$q^{1/2<\deg(\;),\deg(\;)>}$.

We have the following simple lemma
\begin{Lem}\label{bimodV} The algebra $A$ defined in definition \ref{representationring}
is commutative in $U_q-\grmod_P$.
\end{Lem}
\noindent {\it Proof.} The coquasitriangularity of $\Oq$ implies
that $\Oq$ is commutative in the category of $U_q\otimes
U^{\operatorname{op}}_q$-modules (with the obvious braiding). The
subalgebra $A \cong \Oq^{N_q}$ of $\Oq$ is no longer an
$U_q\otimes U^{\operatorname{op}}_q$-module, but an object in
$U_q-\grmod_P$ and the braiding of $U_q\otimes
U^{\operatorname{op}}_q$ acts as the braiding in $U_q-\grmod_P$
making it a commutative algebra there. $\Box$

The $G_q$-equivariant objects in $\Proj(A)$ are by definition
those that corresponds to ${}_{G_q}\MBG$ under the equivalence in
corrolary \ref{proj}. The following result will be useful in the
next section

\begin{Cor}\label{bimodV2}
Any $G_q$-equivariant $M$ in $\Proj(A)$ is an $A$-bimodule.
\end{Cor}
\noindent {\it Proof.} Note that $G_q$-equivariant objects in
$\Proj(A)$ can be thought of as graded $A$-modules with a
compatible $\Oq$-comodule structure. By lemma \ref{bimodV} and the
previous  discussion it follows that they are $A$-bimodules.
$\Box$ This way, we get an action
\begin{equation}\label{actis}
{}_{G_q}\MBG \otimes \MBG \to \MBG,\; M\times N \to
M\otimes_{\Oq}N
\end{equation}
This suggestive notations indicates (ofcourse) that one can define
an $\Oq$-bimodule structure on ${}_{G_q}\MBG$ but we didnt work
this out.

\section{D-modules on Quantum flag variety}
\subsection{Ring of differential operators on $G_q$}

Recall the $U_q$-bimodule structure on $\Oq$ given by \ref{23}.
\begin{defi} We define the ring of quantum differential operators on $G_q$
to be the smash product algebra $\D := \Oq \star U_q$. So $\D =
\Oq \otimes U_q$ as a vector space and multiplication is given by
\begin{equation}\label{e7}
a \otimes u \cdot b \otimes v = au_1(b)\otimes u_2v.
\end{equation}
\end{defi}
We consider now the ring $\D$ as a left $U_q$-module, via the left
$U_q$-action on $\Oq$ in \ref{23} and the left adjoint action on
itself. (This is not the action induced from the ring embedding
$U_q \to 1 \otimes U_q \subset \D$.) This way $\D$ becomes a
module algebra for $U_q$:
\begin{equation}\label{e20}
u \cdot a \otimes v = au_1(b)\otimes u_{21}vS(u_{22}).
\end{equation}
In the following we will use the restriction of this action to
$U_q(\bb)$. As $U_q(\g)$ is not locally finite with respect to the
adjoint action on itself, this $U_q(\bb)$-action doesn't integrate
to a $B_q$-action. Thus $\D$ is not an object of $\MBG$; however,
$\D$ has a subalgebra $\D^{fin} = \Oq \star \U$ which belongs to
$\MBG$. This fact will be used below.
\subsection{$\D$-modules on flag variety}
Let $\lambda \in T_P$.
\begin{defi}\label{d1} A $(B_{q},\lambda)$-equivariant $\D$-module
is a triple $(M, \alpha, \beta)$, where $M$ is a $\C$-module,
$\alpha: \D \otimes M \to M$ a left $\D$-action and $\beta: M \to
M\otimes \Oq(B)$ a right $\Oq(B)$-coaction. The latter action
induces an $U_q(\bb)$-action on $M$ also denoted by $\beta$. These
actions are related as follows:
\smallskip

\noindent $i)$ The $U_q(\bb)$-action on $M \otimes \C_\lambda$
given by $\beta \otimes \lambda$ and by $(\alpha\vert_{U_q(\bb)})
\otimes \operatorname{Id}$ coincide.

\noindent $ii)$ The map $\alpha$ is $U_q(\bb)$-linear with respect
to the $\beta$-action  on $M$ and the action on $\D$ that is given
by \ref{e20} .

\smallskip

These objects form a category denoted $\DBG$. There is the
forgetful functor $\DBG \to \MBG$. Morphisms in $\DBG$ are
morphisms in $\MBG$ that are $\D$-linear.
\end{defi}
We define
\begin{defi}\label{d2} $\Dl$ is the maximal quotient of
$\D$ which is an object of $\DBG$.
\end{defi}
Thus, a simple computation shows that $\Dl \cong \D/ \D I$ where
\begin{equation}\label{relevant}
I = \{E_i, K_i -\lambda(K_i);\, 1 \leq i \leq l\}
\end{equation}
Note that $\D I$ is not a two-sided ideal and hence $\Dl$ is not a
ring. We have
\begin{equation}\label{relevantstruct}
\Dl = \Oq \otimes (U_q/ U_q I) \cong \overline{p}^\star(M_\lambda)
\end{equation}
as a vector space. We define the global section functor $\Gamma:
\DBG \to \M$ to be the global section functor on $\MBG$ composed
with the forgetful functor $\DBG \to \MBG$. Thus $\Gamma =
(\;)^{B_q}$, is the functor of taking $B_q$ invariants (with
respect to the action $\beta$) and we obviously have
\begin{equation}\label{globalD}
 \Gamma = \Hom_{\DBG}(\Dl,\;).
\end{equation}
In particular,  $\Gamma(\Dl) = \End_{\DBG}(\Dl)$ is a ring with
multiplication induced from that in $\D$.
\subsection{} \begin{defi} Let $M_\lambda$ be a Verma module with
highest weight $\lambda$ and put $J_\lambda =
\Ann_{\U}(M_\lambda)$. We define $\Ul = \U/J_\lambda$.
\end{defi} We have
\begin{Prop}\label{globsectD}
There is a ring injection $\Ul \to \Gamma(\Dl)$ which is an
isomorphism for  all $q$ except a finite set of roots of unity
depending on the root data.
\end{Prop}

\noindent ${\mathit{Proof\; of\; proposition\; \ref{globsectD}.}}$
There is the natural surjection $\Ul \to M_\lambda$. It induces a
surjective map
\begin{equation}\label{Inds}
 \overline{p}^\star(\Ul) \to \overline{p}^\star(M_\lambda) = \Dl
\end{equation}
Since $\Ul$ is a $G_q$-module, $\Gamma(\overline{p}^\star(\Ul))$
is isomorphic to $\Ul$. Applying $\Gamma$ to \ref{Inds} we  get
the map
\begin{equation}\label{Inds2}
\Ul \to \Gamma(\Dl)
\end{equation}

\smallskip

\noindent {\it The map \ref{Inds2} is injective when $\lambda =
0$.} Let $\Oloc$ be the localization of $\Oq$ defined by De
Concini and Lyubashenko, \cite{DL}. This is an object in $\DOBG$.
Here comes the structures: As a (right) $\Oq(B)$-comodule
\begin{equation}\label{tyy}
\Oloc = \Oq(N) \otimes \Oq(B)
\end{equation}
Thus, the map $\beta: \Oloc \to \Oloc\otimes \Oq(B)$ is given by
the coproduct of $\Oq(B)$. The map $\alpha: \D \otimes \Oloc \to
\Oloc$ is given as follows: the $\Oq$-module structure on $\Oloc$
is the natural one coming from the localization. The $U_q$-action
on $\Oq$ extends to the localization.

The restriction to $\Uzero$ of the $\Gamma(\Dzero)$-action on
$\Gamma(\Oloc)$ comes from the natural right action of $\Uzero$ on
$\Oloc$. It now follows from \ref{tyy} that $\Gamma(\Oloc)$ is
isomorphic to $M^{\star}_0$ as an $\Uzero$-module. The injectivity
claim now follows since $\Uzero$ by definition acts faithfully on
$M_0$ and hence on $M^*_0$.
\smallskip

\noindent {\it The map \ref{Inds2} is an isomorphism for $\lambda
= 0$.} We define a $\mathbb{Z}$-filtration on $U_q$ by putting
$\deg E_i, \deg F_i = 1$ and $\deg K_i = -1$.  Denote by
$\F_j(Object)$ the $j$'th filtered part of a filtered $Object$;
the associated graded object is denoted by
$$\gr(Object) = \oplus \gr_j(Object) = \oplus \F_j(Object)/F_{j-1}(Object).$$
Intersecting our filtration of $U_q$ with $\U$ we get a filtration
on $\U$ satisfying $\F_j(\U) = 0$ for $j < 0$ and $\dim \F_j(\U) <
\infty$ for all $j$. We get (positive) quotient filtrations on
$M_0$ and $\Uzero$. This way, \ref{Inds} and \ref{Inds2} become
filtered maps. We get maps
\begin{equation}\label{filten}
\F_j(\Uzero) \to \F_j(\Gamma(\Dzero)) = \Gamma(\F_j(\Dzero)) =
\Gamma(\overline{p}^\star(\F_j(M_0)))
\end{equation}
and hence maps
\begin{equation}\label{filten2}
\gr_j(\Uzero) \to \gr_j(\Gamma(\Dzero)) \to \Gamma(\gr_j(\Dzero))
= \Gamma(\overline{p}^\star(\gr_j(M_0)))
\end{equation}
Put
$$\mu_j(q) = \dim_\C \gr_j(\Uzero),\; \nu_j(q) = \dim_\C
\Gamma(\Ind(\gr_j(M_\lambda))).$$ By Kostants theorem (\cite{D},
chapter 8) $\mu_j(1) = \nu_j(1)$. By results of Joseph and Letzter
\cite{JL2} $\mu_j(q)$ is constant for all $q$ except a finite set
of roots of unity. By results of \cite{APW} the $\C$-dimension of
the global sections of the induction of a finite dimensional
$B_q$-module does not depend on $q$. Hence, $\nu_j(q)$ is
independent of $q$. Hence \ref{filten2} is an isomorphism for each
$j$. Hence, \ref{Inds2} with $\lambda = 0$ is an isomorphism by
standard arguments.

\smallskip

\noindent {\it The map \ref{Inds2} is an isomorphism for general
$\lambda$.} We have filtrations on $\Ul, M_\lambda$ etc and get
maps corresponding to \ref{filten2}:
\begin{equation}\label{filten3}
\gr_j(\Ul) \to \Gamma(\gr_j(\Dl))
\end{equation}
By Joseph $\gr(\Ul)$ is independent of $\lambda$. Also,
$\Gamma(\gr_j(\Dl))$ is independent of $\lambda$ since it equals
$\Gamma(\overline{p}^\star(\gr_j(M_\lambda)))$. Under these
identifications the map \ref{filten3} is independent of $\lambda$.
Hence, \ref{filten3} and so \ref{Inds2} are isomorphisms. $\Box$.

\begin{Rem}
\smallskip

\noindent $1)$ Note that the object in $\Proj(A)$ corresponding to
$\overline{p}^\star(\U)$ is $A \otimes \U$. This can be given the
structure of an algebra $A\star \U$. Then one can see that our
$\D$-modules becomes a category of objects in $\Proj(A)$ equipped
with a graded action of this algebra and $\lambda$-compatibility.
This relates our work to the work of Tanisaki, \cite{T}.

\smallskip

\noindent $2)$ Differential operators on the big cell and its
translates of quantum $G/B$ gives the algebras of differential
operators of Joseph, \cite{J}.
\end{Rem}

\subsection{Localization.}
From now on we assume that $q$ is not a root of unity.
\begin{Thm}\label{lockie}
For $\lambda \in \h^*$ regular and dominant, $\Gamma: \DBG \to
\Gamma(\Dl)-\operatorname{mod}$ is an equivalence of categories.
\end{Thm}
Our proof is very similar to Beilinson and Bernsteins proof of
this theorem for classical flag-varieties.

In the following discussion $V$ will denote a finite dimensional
$G_q$-module. It is well known that $V$ admits a filtration
\begin{equation}\label{VVV}
0 \subset V_{0} \subset \ldots \subset V_i  \subset \ldots \subset
V_n = V
\end{equation}
of $B_q$-submodules where $V_i/V_{i-1} \cong \C_{\mu_i}$ and
$\mu_i > \mu_j \implies j > i$. (Thus $\mu_0$ is the highest
weight and $\mu_n$ the lowest weight of $V$.)

\begin{Lem}\label{Filt1}
Let $F \in \MBG$ and consider $V \otimes F$ as a $B_q$-module via
the diagonal action. $(\mathbf{a})$ We have $F^{\dim V} \cong V
\otimes F$ as left $B_q$-modules.   $\mathbf{b})$ The filtration
\ref{VVV} induces a $B_q$-filtration
\begin{equation}\label{VVV2}
0 \subset \ldots \subset V_i \otimes F \subset \ldots \subset V
\otimes F
\end{equation}
We have $V_i \otimes F/ V_{i-1}\otimes F \cong \C_{\mu_i} \otimes
F \cong F \otimes \C_{\mu_i} = F(-\mu_i)$ as $B_q$-modules.
\end{Lem}
{\it Proof of lemma \ref{Filt1}} $(\mathbf{a})$ Let $V^{triv}$ be
$V$ with the trivial $B_q$-action and $V^{triv} \otimes F (\cong
F^{\dim V})$ the $B_q$-module with action on the second factor.
The map
\begin{equation}\label{hallo}
 V\otimes F \to V^{triv} \otimes F;\; v \otimes f \to v_1\otimes
 v_2 f
\end{equation}
is a $B_q$-isomorphism. (Its inverse is $v\otimes f \to v_1\otimes
S(v_2)f$.)
\smallskip

\noindent $(\mathbf{b})$ The only statement that needs a proof is
that $\C_{\mu} \otimes F \cong F \otimes \C_{\mu}$. We construct
an $B_q$-isomorphism $\pi:F \to \C_{\mu} \otimes F \otimes
\C_{-\mu}$ by $\pi(f) = q^{-<\mu,\phi>} 1\otimes f \otimes 1$ for
$f \in F^\phi$. Here $F^\phi$ is the $\phi$-weightspace of $F$.
$\Box$

The filtration \ref{VVV2} induces a projection and an injection
$$p_F: V\otimes F \to F(-\mu_n) \;\,\;\, i_F: F \to V\otimes
F(\mu_0)$$ respectively. (Here $i_F$ is the inclusion $F(-\mu_0)
\cong V_0 \otimes F \to V \otimes F$ tensored by $\C_{-\mu_0}$.)
The isomorphism $V\otimes F \cong F^{\dim V}$ induces an
$\Oq$-module structure on $V\otimes F$ making it an object in
$\MBG$. With this structure the maps $i_F$ and $p_F$ are morphisms
in $\MBG$.

\begin{Rem} Another way to see the $\Oq$-module structure on $V \otimes F$ is to
define it as $\overline{p}^\star(V) \otimes_{\Oq} F$ as we did in
section \ref{commutative}.
\end{Rem}

Assume that $F \in \DBG$. Then each $V_i \otimes F$ becomes an
$U_q$-module, by restricting the $\D$-action on $F$ to an action
of its subalgebra $U_q \cong 1 \otimes U_q$ and using the trivial
$U_q$-action on $V_i$. In this case $p_F$ and $i_F$ are
$U_q$-linear.
\begin{Lem}\label{Filt2} Assume that $F \in \DBG$.
$(\mathbf{a)}$ If $\lambda$ is dominant, then $i_F$ has a
splitting that is $U_q$ and $B_q$-linear. $(\mathbf{b)}$ If
$\lambda$ is regular and dominant, then $p_F$ has a splitting that
is $U_q$ and $B_q$-linear.
\end{Lem}
{\it Proof of lemma \ref{Filt2}.} $\mathbf{a)}$ The center
$\mathcal{Z}$ of $U_q$ acts on $V_i\otimes F(\mu_0)/V_{i-1}\otimes
F(\mu_0)$ by the character $\chi_{-\lambda -\mu_0 +\mu_i}$. But
then $\chi_{-\lambda} \neq \chi_{-\lambda -\mu_0 +\mu_i}$ for
$i\neq 0$. Thus, by Harish-Chandra's theorem the map $i_F$ splits
$U_q$-linearly. The compatibility of the $U_q$ and $B_q$-actions
implies that the splitting map is $B_q$ linear as well.

$(\mathbf{b)}$ The center $\mathcal{Z}$ of $U_q$ acts on
$V_i\otimes F/V_{i-1}\otimes F$ by the character $\chi_{-\lambda
+\mu_i}$. But then $\chi_{-\lambda+\mu_n} \neq \chi_{-\lambda
+\mu_i}$ for $i\neq 0$. Again, this implies that the map $p_F$
splits $U_q$-linearly and hence $B_q$-linearly. $\Box$

\smallskip
\begin{Rem} Exactly as in the classical theory the splittings of $i_F$
and of $p_F$ given by lemma \ref{Filt2} are not $\Oq$-linear.
Since the maps are $B_q$-linear and since lemma \ref{enough2}
shows that the cohomologies $R^j\Gamma$ can be computed by taking
injective resolutions of underlying $B_q$-modules, we see that
they induce splittings on cohomologies.
\end{Rem}

 \noindent {\it Proof of theorem \ref{lockie}}
\smallskip

\noindent $i)$ {\it The functor $\Gamma$ is exact.} Let $F \in
\DBG$. We must prove that $R^j\Gamma(F) = 0$. This will follow if
we can prove that for any noetherian $M \in \MBG$ and injection $M
\hookrightarrow F$ in $\MBG$, the induced maps $a: R^j\Gamma(M)
\to R^j\Gamma(F)$ is the zero map for all $j \geq 0$.

Let $V$ be as in lemma \ref{VVV}. Assume that $\mu_0$ is
sufficiently large for $R^j\Gamma(M(\mu_0))=0$ to hold. We get a
commutative diagram
\begin{equation}\label{basic3}
\begin{matrix}
{R^j\Gamma(M)} & {\overset{i_M}{\to}} &
{R^j\Gamma(V\otimes M(\mu_0))} \\
{{\downarrow}_a} & {} & \downarrow_{} \\
R^j\Gamma(F) & \overset{i_F}{\to} & R^j\Gamma(V\otimes F(\mu_0))
\end{matrix}
\end{equation}
Since $R^j\Gamma(V\otimes M(\mu_0))\cong R^j\Gamma(M^{\dim
V}(\mu_0)) = 0$, the composition $i_F \circ a$ is zero. Since
$i_F$ has a section by lemma \ref{VVV2}, $a$ is zero.

\smallskip

\noindent $ii)$ {\it The functor $\Gamma$ is an equivalence of
categories.} Since we know that $\Gamma$ is exact this follows
from general considerations if we can prove that any $F \in \DBG$
satisfies $\Gamma(F) \neq 0$.

Since $p_F$ splits, $\Gamma(F(-\mu_n))$ is a direct summand in
$\Gamma(V \otimes F) \cong \Gamma(F)^{\dim V}$. If $\mu_n$ is
sufficiently negative we have $\Gamma(F(-\mu_n)) \neq 0$. Hence,
$\Gamma(F) \neq 0$. $\Box$

\end{document}